\documentclass[a4paper, 10pt, conference]{ieeeconf}      

\IEEEoverridecommandlockouts                              
\usepackage{amsmath,amssymb,color,graphicx}
\usepackage{hyperref}
\renewcommand\ge{\geqslant}
\renewcommand\geq{\geqslant}
\renewcommand\le{\leqslant}
\renewcommand\leq{\leqslant}
\newcommand{\I}{[0,1]}

\newcommand{\beeq}{\begin{eqnarray}}
\newcommand{\eeeq}{\end{eqnarray}}

\newcommand{\disp}{\displaystyle}

\newcommand{\co}{\mathop{\rm co}\nolimits}

\newcommand{\abs}[1]{\left\vert#1\mathstrut\right\vert}
\newcommand{\norm}[1]{\left\Vert#1\mathstrut\right\Vert}
\newcommand{\Set}[1]{\left\lbrace#1\mathstrut\right\rbrace}

\newcommand{\st}{\,:\,}
\newcommand{\eps}{\varepsilon}

\renewcommand{\epsilon}{\eps}

\renewcommand{\tilde}{\widetilde}

 \newcommand{\xb}{\bar{x}}

\newtheorem{theorem}{Theorem}[section]
\newtheorem{proposition}{Proposition}

\newcommand{\B}{{\mathbb B}}

\newcommand{\R}{{\mathbb R}}

\title{\LARGE \bf A variant of nonsmooth maximum principle for state constrained problems}

\author{Md. Haider Ali Biswas and M.d.R. de Pinho%
\thanks{This work has been supported by the European Union Seventh Framework Programme [FP7-PEOPLE-2010-ITN] under
grant agreement n°264735-SADCO.
The first author is supported by the grant SFRH/BD/63707/2009,  FCT, Portugal.}
\thanks{Haider Ali Biswas is with
ISR and DEEC, Faculdade de Engenharia da Universidade do Porto,
Rua Dr. Roberto Frias , 4200-465 Porto, Portugal  {\tt\small
dee08022@fe.up.pt}}%
\thanks{Maria do Ros{\'a}rio de Pinho is with
ISR and DEEC, Faculdade de Engenharia da Universidade do Porto,
Rua Dr. Roberto Frias , 4200-465 Porto, Portugal  {\tt\small
mrpinho@fe.up.pt}}%
 }

\begin{document}

\maketitle

\thispagestyle{empty}
\pagestyle{empty}


\begin{abstract}
We derive a variant of the  nonsmooth maximum principle for problems with pure state constraints. The interest of our result resides on the nonsmoothness itself since, when applied to smooth problems, it coincides with known results. Remarkably, in the  normal form, our result has  the special feature of being a sufficient optimality  condition for  linear-convex problems, a feature that the classical Pontryagin maximum principle had whereas the nonsmooth version had not.
This work is distinct to previous work in the literature since, for state constrained problems,  we add the Weierstrass conditions to adjoint inclusions using the joint subdifferentials  with respect to the state and the control. Our  proofs use old  techniques developed in \cite{VinPap82}, while appealing to new results  in \cite{CdP10}.

\end{abstract}

\section{INTRODUCTION}

It is commonly accepted that optimal control appears  with  the  publication
of the seminal book \cite{Pont62} where the statement and proof of the  Pontryagin Maximum Principle played a crucial role
 (we refer the reader to the  survey \cite{PeschPl} for an
interesting historic account of the pioneering results). Since then we have witnessed continuous developments.

Generalization of the classical maximum principle  to problems with nonsmooth data
 appeared in 1970's
as mainly the  result of the  work of Francis Clarke (see \cite{Cla76} and references therein).
The nonsmooth maximum principle, nowadays a well established result, was then extended and refined by a number of authors. One of the first attempts to extend it to cover problems with state constraints came up in
 \cite{VinPap82}.

A special  feature of the classical Pontryagin maximum principle is that it is also a sufficient optimality condition for the normal form of the so called linear-convex problems. Regrettably,
 the nonsmooth version had no such feature.
Nonsmooth necessary optimality conditions in the vein of maximum principles were
  proposed in   \cite{dPV95} overcoming this setback. Regrettably  those necessary conditions
  did not include the Weierstrass condition responsible for the very name Maximum Principle.
More recently the setbacks in  \cite{dPV95} were taken care of in \cite{CdP09}  where a new variant  of the nonsmooth maximum principle is derived
by appealing to \cite{Cla05}. As in  \cite{dPV95}, Lipschitz continuity of dynamics with respect to both
state and control is assumed, the special ingredient responsible for sufficiency of the nonsmooth maximum principle when applied to  normal linear convex problems (see problem $(LC)$ below) \footnote{
With  respect to generalizations of  \cite{dPV95} we also  refer the reader to a different  version of a nonsmooth maximum principle in \cite{Zvi11} making use of  ``compatible'' feedback controls.}.
In what follows, and for simplicity,  we opt to refer to the statement of this
new nonsmooth maximum principle  stated as Theorem 3.1 in \cite{CdP10} which plays a crucial role in our developments.

Here we extend Theorem 3.1 in \cite{CdP10} to cover state constrained problems. In doing so we
follow closely the approach of \cite{DFF02} and \cite{DFF05} where
the main result in  \cite{dPV95} is generalized   to cover state constrained problems in two steps; first the convex case is treated  in \cite{DFF02} using techniques based on  \cite{VinPap82} and then   convexity is removed  in \cite{DFF05}.

In this paper we show that  the proofs in \cite{DFF02} and \cite{DFF05}  adapted easily to allow extension of  Theorem 3.1 in \cite{CdP10} to
  state constrained problems. In this way we obtain a  new variant of
  the nonsmooth maximum principle,  improving  on   \cite{DFF05} by adding  the Weierstrass condition  to the previous conditions while keeping the interesting
feature of being a sufficient condition for normal linear-convex problems.

\section{PRELIMINARIES}
\label{sec:1}

\subsection{Notation}
Here and throughout $\B$ represents the closed unit ball centered
at the origin regardless of the dimension of the underlying space and  $\abs{\,\cdot\,}$ represents the Euclidean norm or the induced matrix
norm on $\R^{p\times q}$.
  The {\em Euclidean
distance function} with respect to a given set $A \subset \R^{k}$
is
$$
d_{A}\colon\R^{k}\rightarrow\R, \qquad
y \mapsto d_{A}(y) = \mbox{inf }\Set{ \abs{y-x} \st x \in A}.
$$

A function $h\colon\I\to\R^p$ lies in $W^{1,1}(\I;\R^p)$
if and only if it is absolutely continuous;
in $L^1(\I;\R^p)$ iff it is integrable; and
in $L^{\infty}(\I;\R^p)$ iff it is essentially bounded. The norm of $L^1(\I;\R^p)$ is denoted by $\norm{\cdot}_1$
and the norm of $L^\infty(\I;\R^p)$ is $\norm{\cdot}_\infty$.



We make use of standard concepts from nonsmooth analysis.
Let $A\subset \R^k$ be a closed set with $\xb\in A$.
The  {\em limiting normal cone to $A$ at $\xb$} is denoted by
$N_{A}(\xb)$.

Given a lower semicontinuous function
$f\colon\R^k\rightarrow\R\cup\Set{+\infty}$
and a point $\xb\in \R^{k}$ where $f(\xb) < +\infty $,
$\partial f(\xb)$ denotes the
{\em limiting subdifferential} of $f$ at $\xb$.
When the function $f$ is Lipschitz continuous near $x$,
the convex hull of the limiting subdifferential, $\mbox{co }\partial f(x)$,
coincides with the {\em (Clarke) subdifferential}.
Properties of Clarke's subdifferentials (upper semi-continuity,
sum rules, etc.), can be found in \cite{Cla83}.
For details on such nonsmooth analysis concepts, see \cite{Cla83}, \cite{rocka:1},  \cite{Vin00}
and \cite{Mor05}.
\bigskip

\subsection{The Problem}
 Consider the problem denoted throughout by $(P)$ of minimizing
 $$l(x(a),x(b))+\displaystyle\int_a^b L(t,x(t),u(t))~dt$$
 subject to the differential equation
 $$\dot{x}(t)= f(t,x(t),u(t))\quad \text{ a.e. }~~t \in [a,b],$$
 the state constraint
 $$h(t,x(t))\leq 0\quad \text{ for all }~~t \in [a,b],$$
 the boundary conditions
 $$(x(a),x(b)) \in C,$$
 and the control constraints
 $$u(t) \in U(t) \quad \text{ a.e. }~~t \in [a,b].$$
Here the interval $[a,b]$ is fixed. We have the state $x(t) \in \R^n$ and the control $u(t)\in \R^k$. The function  describing the dynamics
 is $f:[a,b]\times\R^n\times \R^k\to \R^n$. Moreover $h$ and $ L$  are  scalar functions $h:[a,b]\times\R^n\to \R$,  $L:[a,b]\times\R^n\times \R^k\to \R$, $U$ is a multifunction  and $C\subset
 \R^n\times \R^n$.

We shall denote by $(S)$ the problem one obtains  from $(P)$  in the absence of the  state constraint $h(t,x(t))\leq 0$  and we refer to it as a \emph{standard optimal control problem}.

 Throughout this paper we assume that the following basic assumptions are in force:
 \begin{itemize}
 \item[\bf{B1}]the functions $L$ and $f$  are ${\cal L}\times {\cal B}$-measurable,
\item[\bf{B2}] the multifunction $U$ has ${\cal L}\times {\cal B}$-measurable graph,
\item[\bf{B3}] the set $C$ is closed and $l$ is
locally Lipschitz.
\end{itemize}

For $(P)$ (or $(S)$)  a pair $(x,u)$ comprising an absolutely continuous function $x$, the state, and a measurable function $u$, the control, is called an \emph{admissible process} if it satisfies all the constraints.

An admissible process $(x^*,u^*)$ is a \emph{strong local minimum} of $(P)$ (or $(S)$) if  there exists $\varepsilon >0$ such that
$(x^*,u^*)$ minimizes the cost over all
admissible processes $(x,u)$ such that
\begin{equation}\label{strong}|x(t)-x^*(t)|\leq \varepsilon \text{ for all } t \in [a,b].\end{equation}
It is  a  \emph{local $W^{1,1}$-minimum} if there exists some
 $\epsilon >0$ such that it minimizes the cost to all processes $(x,u)$ satisfying \eqref{strong} and
\beeq\nonumber\displaystyle\int_a^b |\dot{x}(t)-\dot{x}^*(t)|~dt\leq \varepsilon.\nonumber\eeeq
Let $R:[a,b]\to ]0,+\infty]$ be a given measurable function. Then the admissible process $(x^*,u^*)$ is a \emph{ local minimum of radius} $R$
if
it minimizes the cost over all
admissible processes $(x,u)$ such that
\beeq\nonumber& |x(t)-x^*(t)|\leq \epsilon,\qquad |u(t)-u^*(t)|\leq R(t)\quad \hbox{a.e.}\\\nonumber & \text{and}\quad  \displaystyle\int_a^b |\dot{x}(t)-\dot{x}^*(t)|~dt\leq \varepsilon\eeeq
for some $\epsilon>0$.


\subsection{Assumptions}
In what follows the pair $(x^*,u^*)$ will always denote the solution of the optimal control problem under consideration.

Let us take any function $\phi$ defined in $[a,b]\times\R^n\times \R^k$ and taking values in $\R^n$ or $\R$.
\begin{itemize}
\item[\bf{A1}] There exist constants $k_x^\phi$ and $k_u^\phi$ for almost every $t\in [a,b]$ and every $(x_i,u_i)$ ($i=1,2$) such that
$$x_i \in \{x:|x-x^*(t)|\leq  \epsilon\}, \qquad u_i\in U(t)$$ we have
$$|\phi(t,x_1,u_1)-\phi(t,x_2,u_2)|\leq k_x^\phi|x_1-x_2|+k_u^\phi|u_1-u_2|.$$
\item[\bf{A2}] The set valued function $t\to U(t)$ is closed valued  and there exists a constant $c>0$ such that
 for almost every $t\in [a,b]$ we have
$$|u(t)|\leq c\quad \forall u\in U(t).$$
\end{itemize}
When  {\bf A1} is imposed on $f$ and/or  $L$, then the Lipschitz constants are denoted by $k_x^f$, $k_u^f$, $k_x^L$ and $k_u^L$.
Observe that if $U$ is independent of time, then {\bf A2} states that the set $U$ is compact.
Assumption {\bf A2} requires the  controls to
be bounded, a strong hypothesis but nevertheless quite common in applications. It also  simplifies the proofs of the forthcoming results where limits of
sequence of controls needed to be taken.
\vspace{-6pt}

\subsection{Auxiliary Results}

Attention now goes to problem $(S)$, i.e., we assume that the state constraint
is now absent. We next state
 an adaptation of Theorem 3.1 in \cite{CdP10} essential to our analysis in the forthcoming sections.
It is ``an adaptation'' because it holds under stronger assumptions than those in \cite{CdP10}.

\medskip

\begin{theorem}\label{CdP}\emph{
Let $(x^*,u^*)$ be a strong local minimum for problem~$(S)$. If {\bf B1}--{\bf B3} are satisfied,   $f$ and $L$ satisfy {\bf A1} and  $U$ is closed valued, then  there exist
$p\in W^{1,1}([a,b];\R^n)$
  and a scalar
$\lambda_0 \geq 0$ satisfying
the \emph{nontriviality condition} [NT]:
$$||p||_{\infty}+\lambda_0 > 0,$$
the \emph{Euler adjoint inclusion} [EI]:
\beeq & \nonumber
 (-\dot{p}(t),0)\in \partial^C_{x,u}\Big(\langle p, f\rangle-\lambda_0 L\Big)(t,x^*(t),u^*(t))\\
 & \nonumber - \{0\}\times K|p(t)|\partial^C_{u} d_{U(t)}(u^*(t))~~ \text{a.e.},\eeeq
 the global \emph{Weierstrass condition} [W]:\\
\hspace*{1cm} $\forall~~ u\in U(t),$
\beeq & \nonumber
 \langle p(t),f(t,x^*(t),u)\rangle+\lambda_0 L(t,x^*(t),u) \leq \\
 & \nonumber \langle p(t),f(t,x^*(t),u^*(t))\rangle +\lambda_0 L(t,x^*(t),u^*(t))
~~\text{a.e.},\eeeq
 and the \emph{transversality  condition} [T]: $$(p(a), -p(b)) \in  N_{C}^L(x^*(a),x^*(b)) + \lambda_0\partial^L l(x^*(a),x^*(b)).$$}
 In the above $K$ is a constant depending merely on $k_x^f$, $k_x^L$,  $k_u^f$ and  $k_u^L$.
\end{theorem}

In  \cite{CdP10}  the analysis is done for \emph{local minimum of radius} $R$ instead of strong minimum  and it holds under a weaker  assumption than
{\bf A1}.
%
%
%

We point out that  the conditions given by the classical nonsmmoth maximum principle (see \cite{Cla05}) are [NT], [W], [T] and
[EI] is replaced by
\beeq\label{cMP}
& -\dot{p}(t)\in  \\
& \nonumber \partial^C_{x}\Big(\langle p(t), f(t,x^*(t),u^*(t))\rangle-\lambda_0 L(t,x^*(t),u^*(t))\Big).
\eeeq
We refer the reader to \cite{CdP09} for a discussion on \eqref{cMP} and [EI].

\section{MAIN RESULTS }

\vspace{-2pt}

We now turn to problem $(P)$. We derive a new nonsmooth maximum principle for this state constrained problem in the vein of
Theorem 3.1 in \cite{CdP10} in two stages. Firstly the result is established  under a convexity assumption on the ``velocity set'' (see {\bf C} below). Then  such hypothesis is removed. This is proved following an approach in \cite{Vin00}  and similar to what is done in \cite{DFF05}.

\medskip

On $h$ we impose the following:
\begin{itemize}
\item[\bf{A3}] For all $ x$ such that $|x(t)-  x^*(t)|\leq  \epsilon$ the function $t \to h(t, x)$ is continuous. Furthermore,
there exists a constant $k_h>0$  such that the function $x \to h(t, x)$ is
Lipschitz of rank $k_h$ for all $t \in [a,b]$.%

\end{itemize}
The need to  impose continuity of $t\to h$ instead of merely semi upper continuity is discussed in \cite{DFF02}.

Recall that our basic assumptions
{\bf B1}--{\bf B3} are in force. Suppose that $f$ and $L$  satisfy {\bf A1}  and that {\bf A2} holds.
For future use, observe that these assumptions also assert that following conditions are
satisfied:
\beeq\label{Lip}
& |\phi(t,x^*(t),u)-\phi(t,x^*(t),u^*(t))|\leq  \\
\nonumber &  k_u^\phi|u-u^*(t)|\text{ for all }u \in U(t)\text{ a.e. }t
\eeeq
  and there exists an integrable
function $k$ such that
\begin{equation}\label{boundoff} | \phi(t,x^*(t),u)|\le k(t) \text{ for all }u \in U(t)\text{ a.e. }t .
\end{equation}
In the above $\phi$ is to be replaced by  $f$ and  $L$.
Moreover, it is a simple matter to see that the sets
$f(t,x, U(t))$ and $ L(t,x, U(t))$
are compact for all $x\in x^*(t)+\epsilon \B$.

\subsection{Convex Case}

Consider the additional  assumption  on the ``velocity set'':

\begin{itemize}
\item[\bf{C}] The velocity set
$$\left\{ (v,l)=(f(t,x,u),L(t,x,u)),~ u\in U(t)\right\}$$
is convex for all $(t,x)\in [a,b]\times\R^n$.
\end{itemize}

\medskip

Introduce the following subdifferential
\beeq\label{barpartialh}
 & \bar{\partial}_x h(t,x) :=\\  &\nonumber
  \co \{\lim \xi_i\ :\ \xi_i \in \partial_x h(t_i,x_i), (t_i,x_i) \to (t,x) \}.
\eeeq

\medskip

\begin{proposition}\label{DFF}\emph{
Let $(x^*,u^*)$ be a strong local minimum for problem~$(P)$. Assume that  $f$ and $L$ satisfy {\bf A1}, assumptions {\bf B1}--{\bf B3},  {\bf A2} and {\bf C} hold and  $h$ satisfies  {\bf A3}. Then there exist
$p\in W^{1,1}([a,b];\R^n)$,  $\gamma\in L^1([a,b];\R)$,
     a  measure $\mu \in C^\oplus ([a,b];\R)$,
  and a scalar
$\lambda_0 \geq 0$ satisfying
\begin{itemize}
\item[(i)] $\mu \{[a,b]\} + ||p||_{\infty}+\lambda_0 > 0,$\\[0.5mm]
\item[(ii)]
 $(-\dot{p}(t),0)\in $\\
 $\partial^C_{x,u}\Big(\langle q(t),f(t,x^*(t),u^*(t))\rangle-\lambda_0 L(t,x^*(t),u^*(t))\Big)$\\
\hspace{1cm}$-\{0\}\times N_{U(t)}^C(u^*(t)) \text{ a.e.},$\\[0.5mm]
\item[(iii)] $\forall~u\in U(t),$\\
$ \langle q(t),f(t,x^*(t),u)\rangle -\lambda_0 L(t,x^*(t),u) \leq$\\
$ \langle q(t),f(t,x^*(t),u^*(t))\rangle-\lambda_0 L(t,x^*(t),u^*(t))
~\text{a.e.},$\\[1mm]
\item[(iv)] $(p(a), -q(b)) \in $\\
$ N_{C}^L(x^*(a),x^*(b)) + \lambda_0\partial l(x^*(a),x^*(b)),$\\[0.5mm]
\item[(v)] $\gamma(t) \in \bar{\partial} h(t,x^*(t)) \quad \mu \mbox{-} \text{a.e.},$\\[0.5mm]
\item[(vi)] $\mathrm{supp}\{ \mu \} \subset
     \left\{ t \in [a,b] : h(t,x^* (t) )=0  \right\},$
\end{itemize}
where
\begin{equation}
    \label{qqq}
    q(t)=\left\{\begin{array}{l} p(t)+\int_{[a,t)} \gamma( s) \mu (ds)\quad   t \in [a,b)\\[3mm]
                p(t)+\int_{[a,b]} \gamma( s) \mu (ds) \quad   t =b.\end{array}\right.
  \end{equation}}
\end{proposition}

\subsection{Maximum Principle in the Nonconvex Case}
\vspace{-2pt}

Now we replace  the subdifferential $\bar{\partial}_x h$ by a more refined subdifferential $\partial^>_{x}h$ defined by
\beeq\label{partialh}
  & \partial^>_{x} h(t,x) :=\co \{\xi :\exists (t_i,x_i) \xrightarrow{h} (t,x) : \\
  &\nonumber h(t_i,x_i)>0 ~\forall i,  ~\partial_x h(t_i,x_i) \to \xi  \}.
\eeeq

\begin{theorem}\label{DFF1}
\emph{Let $(x^*,u^*)$ be a strong local minimum for problem~$(P)$. Assume that $f$ and $L$ satisfy {\bf A1}, $h$ satisfies {\bf A3} and that
 {\bf A2} as well as the basic assumptions {\bf B1}--{\bf B3} hold.
Then there exist an absolutely continuous function
$p$, an integrable function $\gamma$,
      a non-negative measure $\mu \in C^\oplus([a,b];\mathbb{R})$,
  and a scalar
$\lambda_0 \ge 0$ such that  conditions (i)--(vi) of Proposition \ref{DFF}
hold with $\partial^>_{x}h$ as in (\ref{partialh}) replacing $\bar{\partial}_x h$ and where $q$ is as defined in \eqref{qqq}.}
\end{theorem}

\medskip

For the convex case see \cite{Haider2011} for preliminary results for problems with additional mixed state control constraints. Removal of convexity will the be focus of future work.

The  above theorem adapts easily when we assume  $(x^*,u^*)$ to  be a weak local minimum instead of a strong local minimum (see discussion above). It is sufficient to replace $U(t)$ by $U(t)\cap \B_\epsilon(u^*(t))$.

\medskip

Theorem \ref{DFF1} can now be extended to deal with a  local $W^{1,1}$-minimum for $(P)$.
\begin{theorem}\label{DFF2}\emph{
Let $(x^*,u^*)$ be merely a  local $W^{1,1}$-minimum for problem ~$(P)$. Then the conclusions of Theorem \ref{DFF1} hold.}
\end{theorem}

\medskip

We omit the proof of this Theorem here since it can be easily obtained mimicking what is done in \cite{Vin00}.

\subsection{Linear Convex Problems}

The distinction between Theorem \ref{DFF1} and classical nonsmooth maximum principle (see \cite{Vin00}) is well illustrated by
an example provided in \cite{DFF02}. We recover such example here showing that Theorem \ref{DFF1} can eliminate processes whereas the classical nonsmooth maximum principle cannot.

\medskip

{\bf Example:} Consider the problem on the interval $[0,1]$:
$$
\mathrm{(L)}\ \left\{
 \begin{array}{l}\text{Minimize }  \disp\int_0^1 (w_1 |x- u_1| + w_2|x-u_2| + x) dt  \\
\text{subject to}   \\
\begin{array}{l}
\dot{x}(t)  =    4w_1(t)u_1(t)+4w_2(t)u_2(t) ~\text{ for a. e. } t, \\
x(t)\geq -1 ~\text{ for all }t,\\
u_1(t), u_2(t) \in  [-1, 1] ~\text{for a. e. } t,\\
 (w_1(t),w_2(t))\in W  ~\text{for a. e. } t,\\
x(0) =0
\end{array}
\end{array}
\right.
$$
where
$$W:=\{(w_1,w_2)\in \R^2: ~w_1,~w_2\geq 0,~w_1+w_2=1\}.$$
The process $(x^*,u_1^*,u_2^*,w_1^*,w_2^*):=(0,0,0,1,0)$ is an admissible process with cost $0$
and along the trajectory the state constraint is inactive. It is easy to see that the classical nonsmooth maximum principle holds when we take all the multipliers $0$ but $\lambda_0=1$. However,  $(x^*,u_1^*,u_2^*,w_1^*,w_2^*)$ is not optimal. In fact, if we consider
the process  $(x,u_1,u_2*,w_1,w_2)=(-4\alpha t, -\alpha, 0,1,0)$, with $\alpha \in (0,1/4)$, we see that this process  has cost
$-3/4\alpha$.
Now let us apply Theorem \ref{DFF1} to our problem for the process $(x^*,u_1^*,u_2^*,w_1^*,w_2^*)$.
Since the state constraint is inactive, we deduce that measure $\mu$ is null. Considering
the Euler Lagrange equation in (ii) of Theorem \ref{DFF1} we deduce that there should exists an absolutely
continuous function $p$ and a scalar $\lambda_0\geq 0$ satisfying (i) of Theorem \ref{DFF1} and  such that $p(1)=0$, $-\dot{p}(t)=-\lambda_0(1+e(t))$ and $0=4p(t)+\lambda_0 e(t)$ where $e(t)$ takes values in $[-1,1]$\footnote{The function $e$  appears from the subdifferential of the cost which is clearly nonsmooth due to the presence of the modulus.}. A simple analysis will convince the reader that this situation is  impossible.
This means that Theorem \ref{DFF1} does not hold excluding $(x^*,u_1^*,u_2^*,w_1^*,w_2^*)$ as a minimum.
$~~\blacksquare$
\medskip

Consider the problem
$$
\mathrm{(LC)}\ \left\{
 \begin{array}{l}\text{Minimize } l(x(a),x(b))+ \disp\int_a^b L(t,x(t),u(t)) dt  \\
\text{subject to}   \\
\begin{array}{l}
\dot{x}(t)  =    A(t)x(t)+B(t)u(t) \text{ for a. e. } t\in [a,b], \\
D(t)x(t)\leq 0 ~\text{ for all }t\in  [a, b],\\
u(t) \in U(t) ~\text{for a. e. } t\in  [a, b],\\
(x(a),x(b)) \in E
\end{array}
\end{array}
\right.
$$
where $E$ is convex, the multifunction $U$ is convex valued, the functions $l$ and $(x,u) \to L(t,x,u)$ are convex, the function $A : [0, 1] \to \R^{n\times n}$ is integrable, the function $B : [0, 1] \to \R^{n\times k}$ is measurable, and the function $D : [0, 1] \to \R^{1\times n}$ is continuous. Then $(LC)$  is what we refer to as a linear convex problem with state constraints.

Theorem \ref{DFF1} (and of course Theorem \ref{DFF2}) keeps the significant feature of being a sufficient condition of optimality in the normal form for  problem $(LC)$.
This follows directly from the observation that the proof of  Proposition 4.1 in \cite{DFF02} proves  our claim. No adaptation is required in this case. For completeness we state such proposition here.

We say that a process $(x^*,u^*)$ is a {\bf normal} extremal if it satisfies the conclusions of Theorem \ref{DFF1} with $\lambda_0=1$.

\medskip

\begin{proposition}\label{DFF02} (\cite{DFF02}) If the process $(x^*,u^*)$ is
 a normal extremal for problem $(LC)$, then
it is a minimum.
\end{proposition}

\medskip

Let us return to our previous example.
Problem $(L)$ is what we call a  linear convex problem. It is now obvious that  the process $(x^*,u_1^*,u_2^*,w_1^*,w_2^*):=(0,0,0,1,0)$
does not satisfy the conclusions of Theorem \ref{DFF1}, if it did, then it would be a minimum as asserted by Proposition \ref{DFF02} and it is not.

\section{PROOFS OF THE MAIN RESULTS}

Since our proofs are based on those in \cite{DFF02} and \cite{DFF05} we  we only give a brief sketch
of them, refereing the reader to the appropriate literature for details.

All the results are proved assuming that $L\equiv 0$. The case of $L\neq 0$ is treated by a standard and well known technique.

\medskip

\subsection{Sketch of the Proof of Proposition \ref{DFF}}

\begin{itemize}
\item First the validity of the Proposition is  established  for the simpler problem

$$
\mathrm{(Q)}\hspace{.2in}\left\{
 \begin{array}{l}\text{Minimize }  l(x(b))  \\
\text{subject to}   \\
\begin{array}{l}
\dot{x}(t)  = f(t,x(t),u(t)) ~\text{a.e.}  t\in [a,b] \\
         u(t)  \in  U(t) ~\text{a.e.}  t\in [a,b] \\
     h(t,x(t))  \leq  0       ~ \text{for all }  t\in [a,b] \\
(x(a),x(b))\in  \{x_a\}\times E_b.
\end{array}
\end{array}
\right.
$$
Problem $(Q)$ is a special case of $(P)$ in which
$E=\{x_a\}\times E_b$ and $l(x_a,x_b)=l(x_b)$.

Our proof consists of the following steps
\begin{itemize}
\item[Q1] Define a sequence of problems penalizing the  state-constraint violation. The sequence of problems is
$$(Q_i)
\quad
\left\{
 \begin{array}{l}\text{Minimize }  l(x(b)) + i \displaystyle \int_a^b h^+(t,x(t))~dt \\
\text{subject to}   \\
\begin{array}{l}
 \dot{x}(t)  = f(t,x(t),u(t))\quad  \text{a.e.} ~ t\in [a,b] \\[1mm]
(x(a),x(b)) \in \{x_a\}\times E_b,
\end{array}
\end{array}
\right.
$$
where $
 h^+(t,x):=\max \{ 0, h(t,x)\}.$
 \item[Q2] Assume that
{\bf [IH]} $\disp\lim_{i\to \infty}
\inf\{Q_i\} = \inf\{Q\}.$

 \item[Q3] Set $W$ to be the set of  measurable functions
$u:[a,b]\to\R^k,~~u(t)\in U(t) ~\text{a.e.}$   such that a solution of the differential equation
$
\dot{x}(t)=f(t,x(t),u(t))$, for almost every $t\in [a,b]$, with
   $ x(t) \in   x^*(t)+\varepsilon \B$ for all $t\in [a,b]$ and  $x(a)=x_a$ and $x(b)\in E_b$.
We provide $W$ with the $L^1$ metric defined by $
 \Delta(u,v):=\parallel u-v\parallel_{L_1}
$
and set
$$
J_i(u):=l(x(b)) + i\displaystyle\int_a^b h^+(t,x(t))\,dt.
$$
Then  $(W,\Delta)$
is a complete metric space in which the functional
$J_i\colon W\rightarrow\R$ is continuous.
 \item[Q4] Apply Ekeland's theorem to the sequence of problems of the form
 $$(O_i)~\left\{\begin{array}{rl}
 \text{Minimize} & J_i(u)\\
 \text{subject to } & u\in W
 \end{array}\right.$$
which are closely related to $(Q_i)$.

 The conclusion of application of Ekeland's theorem  shows that
 $(x_i,u_i)$ solves the following  optimal control problem:
$$
(E_i)\quad\left\{
\begin{array}{l}
\text{Minimize }  l(x(b)) + i\displaystyle\int_a^b h^+(t,x(t))\,dt+ \\
\qquad \sqrt{\varepsilon_i}\int_a^b \abs{u(t)-u_i(t)}~ dt \\
\text{subject to}   \\
\begin{array}{l}
 \dot{x}(t) = f(t,x(t),u(t))    ~ \text{a.e. }  t\in [a,b] \\
      u(t)\in  U(t)        ~ \text{a.e. }  t\in [a,b] \\
 x(a)=  x_a     \\
   x(b)\in  E_b.
\end{array}
\end{array}
\right.
$$
The fact  that $\varepsilon_i\to 0$ allows us to prove that $u_i$ converges strongly to $u^*$ and $x_i$ converges uniformly to
$x^*$.

 \item[Q6] Rewriting these conditions and taking limits as in \cite{DFF02} we get the required conclusions. 
\item[Q7] Finally we show that {\bf C} implies IH.
\end{itemize}
The remaining of the proof has  three stages.
We first extend Proposition \ref{DFF} to problems where
$x(a) \in E_a,$ and $E_a$ is a closed set.
This is done following the lines  in the end of the proof of Theorem 3.1 in \cite{VinPap82}.

Next we consider the case when the cost is $l=l(x(a),x(b))$. This is done using the technique
in  Step 2 of section 6 in \cite{DLS09}.
And finally,  following again the approach in  section 6 in \cite{DLS09}, we
derive necessary conditions when
$(x(a),x(b)) \in E$, completing the proof.
\end{itemize}

In order to proof our result, an  important piece of analysis added to the proof of Theorem 3.1 in \cite{DFF02} concerns the Weierstrass condition (iii) of Proposition \ref{DFF}.
The information extracted while taking limits allow us to do that  without that much ado.

\subsection{Sketch of the Proof of Theorem \ref{DFF1}}
We now proceed to prove our main Theorem \ref{DFF1}. We recall  that under our  hypotheses both \eqref{Lip}
and \eqref{boundoff} hold and that the set $f(t,x,U(t))$ is compact.

Our proof consists of several steps. We first consider the following 'minimax' optimal
control problem where the state constraint functional $\disp\max_{t\in[a,b]} h(t, x(t))$ appears in the cost.

$$(\tilde{R})~
\left\{
 \begin{array}{l}\text{Minimize }  \tilde{l}(x(a),x(b),\disp\max_{t\in[a,b]} h(t,x(t)))  \\
\text{over $x\in W^{1,1}$ and measurable  $u$ satisfying }    \\
\begin{array}{l}
\dot{x}(t) =  f(t,x(t),u(t)) ~ \text{a.e.} ~ t\in [a,b] \\
u(t) \in U(t)      ~ \text{a.e.} ~ t\in [a,b] \\
(x(a),x(b)) \in  E_a\times \R^n.
\end{array}
\end{array}
\right.
$$
where $\tilde{l}:\R^n \times \R^n\times \R \to \R$ is a given function and  $E_a \subset \R^n $ is a given closed set. We observe that $(\tilde{R})$ is the optimal control problem with free endpoint constraints.

We impose here the following additional assumption $\bf{A4}$, the necessity of which for the forthcoming development of our proof will  become clear soon.
\begin{itemize}
\item[\bf{A4}] The integrable function $\tilde{l}$ is Lipschitz continuous on a neighbourhood of
$$(x^*(a),x^*(b),\max_{t\in[a,b]} h(t,x^*(t)))$$
and $\tilde{l}$ is monotone in the $z$ variable, in the sense that $ z'\ge z$ implies
$\tilde{l}(y,x,z') \ge \tilde{l}(y,x,z), $
for all $ (y,x) \in \R^n \times \R^n.$
\end{itemize}

\smallskip

The following proposition is a straightforward adaptation of Proposition 9.5.4 of \cite{Vin00}.

\medskip

\begin{proposition}\label{FEC1}
\emph{Let $(x^*,u^*)$ be a strong local minimum for problem ~$(\tilde{R})$. Assume
the basic hypotheses, {\bf A1}, {\bf A2} and {\bf A3}  and the data for the problem $(\tilde{R})$ satisfies the
hypothesis $\bf{A4}$.
Then there exist an absolutely continuous function
$p:[a,b]\to\R^n$, an integrable function  $\gamma:[a,b] \to \R^n$,
      a non-negative measure $\mu \in C^\oplus([a,b];\mathbb{R})$,
  and a scalar
$\lambda_0 \ge 0$ such that
\beeq \label{p:t21}
 & \mu \{[a,b]\} + ||p||_{\infty}+\lambda_0 > 0,\\[2mm]
\label{p:t22}
 & (-\dot{p}(t),0)\in \partial^C_{x,u}\langle q(t), f(t,x^*(t),u^*(t))\rangle \\[2mm]
\nonumber
& -\{0\}\times  N_{U(t)}^C(u^*(t)) ~~\text{a.e. } \\[2mm]
\label{p:t24}
& (p(a), -q(b),\disp\int_{[a,b)} \mu (ds))  \in \\
\nonumber &  N_{C_a}^L(x^*(a)) \times \{0,0 \} + \\
\nonumber &\lambda_0\partial \tilde{l}(x^*(a),x^*(b),\max_{t\in[a,b]} h(t,x^*(t)),\\[2mm]
\label{p:t25}
& \gamma(t) \in \bar{\partial} h(t,x^*(t)) \quad \mu \mbox{-} \text{a.e.},\\[2mm]
\label{p:t26}
& \forall~~ u\in U(t),\\
\nonumber & \langle q(t),f(t,x^*(t),u)\rangle \leq \langle q(t),f(t,x^*(t),u^*(t))\rangle ~~\text{a.e. },\\[2mm]
\label{p:t27}
& \mathrm{supp}\{ \mu \} \subset\\
\nonumber &
     \left\{ t \in [a,b] : h(t, x^*(t))=\max_{s\in[a,b]} h(s,x^*(s))  \right\},
\eeeq}
where $q$ is defined as in (\ref{qqq}).

 \end{proposition}

\medskip

We now turn to the  derivation of Theorem \ref{DFF1}.
Consider the set
\begin{equation}
\begin{array}{l}
V:=\{(x,u,e): (x,u) ~\text{satisfies} ~ \dot{x}(t)=f(t,x(t),u(t)), \\
 ~~ u(t) \in U(t) ~\text{a.e.}, ~e \in \R^n, (x(a),e)\in C ~\\
 \text{and} ~\|x-x^*\|_{L^\infty}\le \epsilon\}
\end{array}
\end{equation}
and let  $d_V: V\times V \to \R$ be a function defined by
\beeq\label{CMS}
& d_V((x,u,e),(x',u',e'))=\\
& \nonumber |x(a)-x'(a)|+ |e-e'|+\int^b_{a}|u(t)-u'(t)|dt
\eeeq
For all $i$, we choose $\epsilon_i \downarrow 0$ and define the function
\beeq \nonumber & \tilde{l}_i(x,y,x',y',z):=\\
 \nonumber &\max \{l(x,y)-l(x^*(a),x^*(b))+\epsilon_i^2, z, |x'-y'|\}.\eeeq
Then  $d_V$ defines a metric on the set $V$ and $(V,d_V)$ is a complete metric space
such that
\begin{itemize}
\item If $(x_i,u_i,e_i) \to (x,u,e)$ in the metric space $(V,d_V)$, then $\|x_i-x\|_{L^\infty} \to 0$,
\item  The function
$$(x,u,e) \to \tilde{l}_i(x(a), e, x(b), e, \max_{t\in [a,b]} h(t,x(t)))$$ is continuous on $(V,d_V).$
\end{itemize}

We now consider the following optimization problem
\beeq \nonumber &\text{Minimize} ~\{\tilde{l}_i(x(a), e, x(b), e, \max_{t\in [a,b]} h(t,x(t))):\\
\nonumber & (x,u,e) \in V\}.
\eeeq

We observe that $$\tilde{l}_i(x^*(a), x^*(b), x^*(b), x^*(b), \max_{t\in [a,b]} h(t,x^*(t)))=\epsilon_i^2.$$

Since  $\tilde{l}_i$ is non-negative valued, it follows that $(x^*, u^*, x^*(b))$ is an $\epsilon_i^2 \text{-minimizer}$
for the above minimization problem.
According to Ekeland's Theorem there exists a sequence $\{(x_i, u_i, e_i)\}$ in $V$ such that for each $i$, we have
\beeq\label{mp}
& \tilde{l}_i(x_i(a), e_i, x_i(b), e_i, \max_{t\in [a,b]} h(t,x_i(t))) \le \\
& \nonumber \tilde{l}_i(x(a), e, x(b), e, \max_{t\in [a,b]} h(t,x(t)))+\\
& \nonumber \varepsilon_id_V((x,u,e),(x_i,u_i,e_i))
\eeeq
for all $(x,u,e)\in V$ and we also have
\begin{equation}\label{fc}
d_V((x_i,u_i,e_i),(x^*,u^*,x^*(b))) \le \varepsilon_i.
\end{equation}

Thus the condition (\ref{fc}) implies that $e_i \to x^*(b) ~\text{and}~ u_i \to u^*$ in the $L^1$ norm. By using subsequence extraction,
we conclude  that $u_i \to u^* ~\text{a.e. and}~ x_i \to x^* ~ \text{uniformly}.$

Now we define the arc $y_i\equiv e_i.$ Accordingly we get $y_i \to x^*(b) ~ \text{uniformly}.$
From the minimization property (\ref{mp}), we say that $(x_i,y_i,w_i\equiv 0,u_i)$ is a strong local minimum for the optimal control problem
$$(\tilde{R_i})~
\left\{
 \begin{array}{l}\text{Minimize }  \\
 \tilde{l}_i(x(a), y(a), x(b), y(b), \max_{t\in [a,b]} h(t,x(t)))  \\
                         ~       +\varepsilon_i[|x(a)-x_i(a)|+ |y(a)-y_i(a)|+w(b)] \\
\text{over $x, y, w\in W^{1,1}$ and measurable  $u$ satisfying }    \\
\begin{array}{l}
\dot{x}(t) =  f(t,x(t),u(t)), \dot{y}(t) =0, \\
\dot{w}(t) = |u(t)-u_i(t)| ~ \text{a.e.},  \\
u(t) \in U(t)      ~ \text{a.e.},  \\
(x(a),y(a),w(a)) \in  C\times \{0\}.
\end{array}
\end{array}
\right.
$$

Now  we observe that the cost function of $(\tilde{R_i})$ satisfies all the assumptions of the Proposition \ref{FEC1} and thus this is an example of optimal control problem where the special case of maximum principle of Proposition \ref{FEC1} applies. Rewriting the conclusions of Proposition \ref{FEC1} and taking limits
we obtained the required conditions. The remain of the proof follows closely the approach in \cite{DFF05}.

\section{CONCLUSIONS}
In this work we derive a variant nonsmooth maximum
principle for state constrained problems. The novelty of this work is that our results are also sufficient conditions of optimality for the normal linear-convex problems.
The result presented in the main theorem is quite distinct to previous work in the literature since for state constrained problems,  we add the Weierstrass conditions to adjoint inclusions using the joint subdifferentials  with respect to the state and the control. The illustrated example presented in the paper justifies our results.

\end{document}